# COMBINATORICS & TOPOLOGY OF STRATIFICATIONS OF THE SPACE OF MONIC POLYNOMIALS WITH REAL COEFFICIENTS

BORIS SHAPIRO AND VOLKMAR WELKER[1]

ABSTRACT. We study the stratification of the space of monic polynomials with real coefficients according to the number and multiplicities of real zeros. In the first part, for each of these strata we provide a purely combinatorial chain complex calculating (co)homology of its one-point compactification and describe the homotopy type by order complexes of a class of posets of compositions. In the second part, we determine the homotopy type of the one-point compactification of the space of monic polynomials of fixed degree which have only real roots (i.e., hyperbolic polynomials) and at least one root is of multiplicity $k$. More generally, we describe the homotopy type of the one-point compactification of strata in the boundary of the set of hyperbolic polynomials, that are defined via certain restrictions on root multiplicities, by order complexes of posets of compositions. In general, the methods are combinatorial and the topological problems are mostly reduced to the study of partially ordered sets.

[1]Supported by DFG through "Habilitationsstipendium" 1479/4, work on the project was completed while the author was "General Member" of the MSRI 1996/97 program on combinatorics, research at MSRI is supported in part by NSF grant DMS-9022140

1991 *Mathematics Subject Classification.* Primary 06A07, 58C27. Secondary 57R45.

1. INTRODUCTION

Let $f(X) = X^n + c_{n-1}X^{n-1} + \cdots + c_0$ be a monic polynomial of degree $n$ with coefficients in the real numbers. The space of all monic polynomials of degree $n$ is homeomorphic to real $n$-space $\mathbf{R}^n$. If we speak of $\mathbf{R}^n$ as the space of monic real polynomials of fixed degree $n$ we denote it by $\mathcal{P}ol^n$. We call $f(X)$ a *hyperbolic* polynomial if all roots of $f(X)$ are real. We call $f(X)$ an *elliptic* polynomial if all roots of $f(X)$ are complex (with non-vanishing imaginary part). In particular, an elliptic polynomial is of even degree. The terminology is borrowed from the theory of partial differential equations, see Hörmander's book [18], where elliptic and hyperbolic polynomials in several variables are defined.

In this paper we are concerned with topological properties of strata from the stratification described below. We study the homotopy type and the homology of these spaces by combinatorial techniques. Some of these questions are reduced to the study of certain subcomplexes of the boundary complex of a simplex that are indexed by number-partitions.

1.1. **Stratification.** Let us describe a cellular decomposition of $\mathcal{P}ol^n$ whose cells are indexed by the compositions of numbers $\leq n$. A *composition* $\pi = (a_1, \ldots, a_t)$ of a number $l$ is an ordered tuple of positive integers such that $a_1 + \cdots + a_t = l$ (i.e., a representation of $l$ as an ordered sum). Let $(a_1, \ldots, a_t)$ be a composition of $l$. We denote by $C_l$ the set of all compositions of the number $l$. We order the set $C_l$ by the partial order "$\preceq$" given by refinement (i.e., the minimal compositions greater than the composition $(a_1, \ldots, a_t)$ are of the form $(a_1, \ldots, a_{i-2}, a_{i-1} + a_i, a_{i+1}, \ldots, a_t)$ for some $2 \leq i \leq t$). If $(a_1, \ldots, a_t)$ is a composition of $l$ then we denote by $\mathcal{H}yp^n_{(a_1,\ldots,a_t)}$ the set of hyperbolic polynomials $f(X) = (X-x_1)^{a_1}\cdots(X-x_t)^{a_t}$ satisfying the strict inequalities $x_1 < \cdots < x_t$. By $\mathcal{H}yp^l$ we denote the space of hyperbolic polynomials in $\mathcal{P}ol^l$. It is a simple observation that $\mathcal{H}yp^l$ is homeomorphic to $\mathbf{R} \times \mathbf{R}_+^{l-1}$, where $\mathbf{R}_+$ denotes the space of real numbers $\geq 0$. Thus $\mathcal{H}yp^l_{(a_1,\ldots,a_t)}$ is an open $t$-dimensional cell in the closed subspace $\mathcal{H}yp^l$ of $\mathbf{R}^l$. Let us denote for an even number $n$ by $\mathcal{E}ll^n$ the space of monic elliptic polynomials of degree $n$. Again, standard reasoning implies that $\mathcal{E}ll^n$ is homeomorphic to $\mathbf{C}^{\frac{n}{2}}$. For an arbitrary number $n$ and a composition $\pi = (a_1, \ldots, a_t)$ of some number $l \leq n$ with $n - l$ even, – if $n$ is even we allow $l = 0$ and $\pi = ()$ – we denote by $\mathcal{P}ol^n_{(a_1,\ldots,a_t)}$ the cell $\mathcal{H}yp^l_{(a_1,\ldots,a_t)} \times \mathcal{E}ll^{n-l}$. We call each $\mathcal{P}ol^n_{(a_1,\ldots,a_t)}$ a *standard cell* in $\mathcal{P}ol^n$. For $n$ even and $\pi = ()$ we obtain $\mathcal{P}ol^n_\pi = \mathcal{E}ll^n$. Now if $\pi$ runs over all compositions of numbers $l \leq n$ with $n - l$ even then the cells $\mathcal{P}ol^n_\pi$ stratify $\mathcal{P}ol^n$.

In a next step we group the compositions of the number $n$ into classes. A *number-partition* $\lambda = (\lambda_1 \leq \cdots \leq \lambda_t)$ of the number $n$ is an ascending sequence of positive numbers such $n = \lambda_1 + \cdots + \lambda_t$. We use the notation $\lambda = (1^{e_1}, \ldots, n^{e_n})$ to denote a number-partition consisting of $e_i$ copies of $i$, $1 \leq i \leq n$. One may regard a number-partition as an equivalence class of compositions. A composition $(a_1, \ldots, a_t)$ is of *type* $\lambda = (1^{e_1}, \ldots, n^{e_n})$ if there are exactly $e_i$ numbers $i$ among the numbers $(a_1, \ldots, a_t)$ for $1 \leq i \leq n$. For a number-partition $\lambda$ of $n$, we denote by $\mathcal{H}yp_\lambda^n$ the closure of the union of all cells $\mathcal{H}yp_{(a_1,\ldots,a_t)}^n$, for $(a_1, \ldots, a_t)$ of type $\lambda$, in $\mathcal{H}yp^n$. For a number-partition $\lambda$ of a number $l \leq n$ such that $n-l$ is even, we write $\mathcal{P}ol_\lambda^n$ for $\mathcal{H}yp_\lambda^l \times \mathcal{E}ll^{n-l}$. We call each $\mathcal{P}ol_\lambda^n$ a *standard stratum* in $\mathcal{P}ol^n$. In particular, if $n$ is even and $\lambda = ()$ then $\mathcal{P}ol_\lambda^n = \mathcal{E}ll^n$.

The main results of the paper are :

- For a number-partition $\lambda = (\lambda_1 \leq \cdots \leq \lambda_t)$ of a number $l$ with $n - l$ even, we describe the homotopy type of the one-point compactification of $\mathcal{P}ol_\lambda^n$ as a double suspension over another space whose topology is "essentially" described via the partially ordered set of compositions of numbers $l \leq n$ with $n-l$ even and whose type is coarser than $\lambda$. In particular, the homotopy type of $\mathcal{P}ol_\lambda^n$ only depends on $n$ and the combinatorics of $\lambda$. Then we present a chain complex calculating its cohomology.
- For a number-partition $\lambda$ of $n$ we describe the one-point compactification $\widehat{\mathcal{H}yp_\lambda^n}$ as a double suspension over the order complex of some partially ordered set of compositions. In case the number-partition $\lambda$ is free of resonances (see Section 3.1) we also describe $\widehat{\mathcal{H}yp_\lambda^n}$ as a quotient of the permutahedron and give the homotopy type.
- The one-point compactification $\widehat{\mathcal{H}yp_{(k,1^{n-k})}^n}$ of the space $\mathcal{H}yp_{(k,1^{n-k})}^n$ of all hyperbolic polynomials with at least one $k$fold root is homotopic to a sphere or contractible.

Our attempts were motivated by results of Arnold [1], [2] and Vassiliev [25] on similarly defined spaces of real and complex polynomials. Spaces of rational functions with real coefficients appear in the realm of control theory. We refer the reader to the work of Brockett [8], the more recent results by Helmke [16], [17] and the references contained therein. The general flavor of our definitions is also influenced by the work of Björner [3] on arrangements of linear subspaces in $\mathbf{R}^n$ and in particular by the class of "orbit arrangements" defined in [3]. The strata $\mathcal{H}yp_\lambda^n$ and $\mathcal{P}ol_\lambda^n$ may be regarded as "unordered analogs" of the arrangements introduced and studied by Björner.

Convexity properties of $\mathcal{H}yp^n$ – in the coefficients of the polynomials as coordinates – have been studied rather intensively. We refer the reader to the papers by Meguerditchian [20] and Dedieu [11] and the references contained therein.

This work war started while the second author enjoyed the hospitality of KTH in Stockholm and finished while the first author was visiting the Max Planck Institute in Bonn. We thank both institutions and A. Björner at KTH for providing these opportunities. Finally, we are grateful to U. Helmke for his interest in our work and for providing helpful references.

## 2. Cohomology of 1-point compactification of $\mathcal{P}ol_\lambda^n$

2.1. **Adjacency of cells and strata.** Let us describe combinatorially the closure of the standard cells $\mathcal{P}ol_{(a_1,\ldots,a_t)}^n$. In order to pursue this program we introduce a partial order on the set $C_{\leq l}$ of all compositions of numbers $m \leq l$ with $l - m$ even extending the partial order "$\preceq$" on the compositions $C_m \subseteq C_{\leq l}$ for $m \leq l$. Let us define the cover relation in $C_{\leq l}$. A composition $(a_1,\ldots,a_t)$ covers a composition $(b_1,\ldots,b_s)$ if and only if either

- $s = t + 1$ and $(b_1,\ldots,b_s)$ arises from $(a_1,\ldots,a_t)$ by adding a block "2" (i.e., there is an $1 \leq i \leq s$ such that $a_j = b_j$, $1 \leq j \leq i - 1$, $b_i = 2$, $a_{j-1} = b_j$, $i + 1 < j \leq s$) or
- $s = t - 1$ and $(b_1,\ldots,b_s)$ arises from $(a_1,\ldots,a_t)$ by merging two neighboring blocks (i.e., $(a_1,\ldots,a_t) \prec (b_1,\ldots,b_s)$ in $C_m$, $m = \sum_{j=1}^{t} a_j = \sum_{j=1}^{s} b_j$).

The partial order on $C_{\leq l}$ is the transitive closure of the above defined cover relation. Since this partial order coincides with the "usual" partial order "$\preceq$" on each $C_m$, $m \leq l$, we also denote this partial order by "$\preceq$." For a number-partition $\lambda = (\lambda_1 \leq \cdots \leq \lambda_t)$ of $m \leq l$ with $l - m$ even we denote by $C_{\lambda, \leq l}$ all compositions in $C_{\leq l}$ that are coarser than some composition of type $\lambda$ in the order "$\preceq$."

Obviously, for any number-partition $\lambda$ of $m \leq n$ with $n - m$ even one has

$$\mathcal{P}ol_\lambda^n = \bigcup_{(a_1,\ldots,a_t)\preceq(b_1,\ldots,b_s)} \mathcal{P}ol_{(b_1,\ldots,b_s)}^n,$$

where the union is over $C_{\lambda, \leq n}$.

**Proposition 2.1.** *Let $(a_1,\ldots,a_t)$, $(b_1,\ldots,b_s) \in C_{\leq l}$ be compositions. The cell $\mathcal{P}ol_{(b_1,\ldots,b_s)}^n$ lies in the boundary of $\mathcal{P}ol_{(a_1,\ldots,a_t)}^n$ if and only if the corresponding compositions satisfy $(b_1,\ldots,b_s) \preceq (a_1,\ldots,a_t)$ in $C_{\leq l}$.*

**Proof.** This is obvious. □

2.2. **Homotopy Type.** The homotopy type of the one-point compactification $\widehat{\mathcal{P}ol_\lambda^n}$ is described via a double suspension.

Consider the obvious action of the *"positive" affine group* $Aff_+(\mathbf{R}) = \{X \mapsto \rho X + \gamma \mid \rho > 0, \gamma \in \mathbf{R}\}$ on the space of monic polynomials of even degree $n$ by

$$X^n + c_{n-1}X^{n-1} + \cdots + c_0 \mapsto \frac{\left((\rho X + \gamma)^n + c_{n-1}(\rho X + \gamma)^{n-1} + \cdots + c_0\right)}{\rho^n}.$$

Let $\Sigma_n$ denote the one-parameter family of polynomials of the form $(X + \alpha)^n$, $\alpha \in \mathbf{R}$.

**Proposition 2.2.** *The action of $Aff_+(\mathbf{R})$ on the space $\mathcal{P}ol^n \setminus \Sigma_n$ is free. The complement $\mathcal{P}ol^n \setminus \Sigma_n$ is homeomorphic to $\mathbf{R} \times \mathbf{R}^+ \times S^{n-2}$, where $S^{n-2}$ is the unit sphere in the space $\widetilde{\mathcal{P}ol}^{n-1}$ of polynomials of the form $X^n + c_{n-2}X^{n-2} + \cdots + c_0$.*

**Proof.** Assume a polynomial $f(X) = X^n + c_{n-1}X^{n-1} + \cdots + c_0$ is stabilized by the map $X \to \rho X + \gamma$ under the action of $Aff_+(\mathbf{R})$. Then

$$c_{n-1} = \frac{c_{n-1} + \gamma}{\rho} \Leftrightarrow c_{n-1}(\rho - 1) = \gamma.$$

Thus if $c_{n-1} = 0$ then $\gamma = 0$. In this case $c_i = \frac{c_i}{\rho^{n-i}}$, $0 \leq i \leq n-2$. But this implies either $\rho = 1$ or $c_i = 0$, $0 \leq i \leq n-2$. Thus either $f(X) \in \Sigma_n$ or $\rho X + \gamma = X$. If $c_{n-1} \neq 0$ then the image of $f(X)$ under the transformation $X \mapsto X - \frac{c_{n-1}}{n}$ has coefficient $0$ at the term $X^{n-1}$. Thus the stabilizer of $f(X - c_{n-1})$ is conjugate to the stabilizer of $f(X)$. In particular, either $f(X) = (X + \alpha)^n$ or $\rho X + \gamma = X$. For each element $f(X) \in \mathcal{P}ol^n \setminus \Sigma_n$ there is a unique element $\rho X + \gamma \in Aff_+(\mathbf{R})$ such that $f(X)^{\rho X + \gamma} = X^n + d_{n-1}X^{n-1} + \cdots + d_0$ satisfies $d_{n-1} = 0$ and $d_{n-2}^2 + \cdots + d_0^2 = 1$. To see this one realizes that $d_{n-1} = 0$ implies $\gamma = c_{n-1}$. Thus we may assume $\gamma = c_{n-1} = 0$. Then $f(X)^{\rho X} = X^n + \frac{c_{n-2}}{\rho^2}X^{n-2} + \cdots + \frac{c_0}{\rho^n}$. For $\rho > 0$ the function $h(\rho) = (\frac{c_{n-2}}{\rho^2})^2 + \cdots + (\frac{c_0}{\rho^n})^2$ is strictly decreasing. Since, $\lim_{\rho \to 0} h(\rho) = \infty$ and $\lim_{\rho \to \infty} h(\rho) = 0$ there is a unique $\rho$ such that $h(\rho) = 1$. Thus, $\mathcal{P}ol^n \setminus \Sigma_n$ is a fibered over $S^{n-2}$ with fiber $\mathbf{R} \times \mathbf{R}_{>0}$. Since the fibration is trivial – the map $(f(X), \rho X + \gamma) \mapsto f(X)^{\rho X + \gamma}$ defines an homeomorphism $S^{n-2} \times Aff_+(\mathbf{R}) \cong \mathcal{P}ol^n \setminus \Sigma_n$ – the assertion follows. □

**Corollary 2.3.**
(a) *For every number-partition $\lambda = (\lambda_1 \leq \cdots \leq \lambda_t)$ of a number $l \leq n$ with $n - l$ even the space $\mathcal{P}ol_\lambda^n \setminus \Sigma_n$ is homeomorphic to $\mathbf{R} \times \mathbf{R}_{>0} \times (\mathcal{P}ol_\lambda^n \cap S^{n-2})$, where $S^{n-2}$ is the unit sphere in the space $\widetilde{\mathcal{P}ol}^n$ of monic polynomials of degree $n$ with vanishing coefficient at $X^{n-1}$. In particular, $\mathcal{P}ol_\lambda^n \setminus \Sigma_n$ and $\mathcal{P}ol_\lambda^n \cap S^{n-2}$ are homotopically equivalent.*
(b) *$\mathcal{P}ol_\lambda^n$ is homeomorphic to $\mathbf{R} \times Cone(\mathcal{P}ol_\lambda^n \cap S^{n-2})$, where $Cone(\mathcal{P}ol_\lambda^n \cap S^{n-2})$ is a cone over $\mathcal{P}ol_\lambda^n \cap S^{n-2}$.*

**Proof.**
(a) This follows immediately from Proposition 2.2. Note that the fibration constructed in Proposition 2.2 respects the standard strata and cells – outside $\overline{\Sigma_n}$.
(b) If $f(X)$ is a monic polynomial of degree $n$. Then $\lim_{\rho \to +\infty} f(X)^{\rho X} = X^n$. Thus each $\mathcal{P}ol_\lambda^n$ is homeomorphic to a cone over $\mathcal{P}ol_\lambda^n \cap S^{n-2}$ with apex $X^n$. □

In the formulation of the next corollary we denote by $Susp(X)$ the suspension of the space $X$.

**Corollary 2.4.** *For every number-partition $\lambda = (\lambda_1 \leq \cdots, \lambda_t)$ of a number $l \leq n$ the one-point compactification $\widehat{\mathcal{P}ol_\lambda^n}$ is the one-point compactification (i.e., the Thom space) of the trivial bundle $\mathbf{R} \times \Sigma(\mathcal{P}ol_\lambda^n \cap S^{n-2})$. In particular, $\widehat{\mathcal{P}ol_\lambda^n}$ is homotopic to the wedge*

$$S^1 \vee Susp(Susp(Ell^{n,\lambda} \cap S^{n-2})).$$

**Proof.** By Corollary 2.3 the space $\mathcal{P}ol_\lambda^n$ is homeomorphic to $\mathbf{R} \times Cone(\mathcal{P}ol_\lambda^n \cap S^{n-2})$. Therefore, $\widehat{\mathcal{P}ol_\lambda^n}$ can be obtained by a sequence of two one-point compactifications $\left(\mathbf{R} \times Cone(\mathcal{P}ol_\lambda^n \cap S^{n-2})\widehat{\phantom{x}}\right)$. The one-point compactification of $Cone(\mathcal{P}ol_\lambda^n \cap S^{n-2})$ is $Susp(\mathcal{P}ol_\lambda^n \cap S^{n-2})$, since $\mathcal{P}ol_\lambda^n \cap S^{n-2}$ is closed and compact. Thus $\widehat{\mathcal{P}ol_\lambda^n}$ is homeomorphic to $\left(\mathbf{R} \times Susp(\mathcal{P}ol_\lambda^n \cap S^{n-2})\right)\widehat{\phantom{x}}$. The one-point compactification of a cylinder over a compact (non-empty) set $X$ is homotopic to a wedge of a circle and the suspension over $X$. This implies the rest of the assertion. □

2.3. **(Co)chain complex.** Let us now describe combinatorially a (co)-chain complex for the calculation of the (co)-homology of $\mathcal{P}ol_\lambda^n$ – or, equivalently by Corollary 2.4 modulo a twofold shift the (co)-homology of $\mathcal{P}ol_\lambda^n \cap S^{n-2}$.

Let us introduce the following orientation on $\mathcal{P}ol_{(a_1,\ldots,a_t)}^n$ given by the top differential form

$$\xi = du_1 \wedge dv_1 \wedge \cdots \wedge du_{(n-\sum a_i)} \wedge dv_{(n-\sum a_i)} \wedge dx_1 \wedge \cdots \wedge dx_t.$$

The notation is chosen in the following way. Let $F(X) = g(X)h(X)$ be a polynomial in $\mathcal{P}ol_{(a_1,\ldots,a_t)}^n$, where $g(X) \in \mathcal{H}yp_{(a_1,\ldots,a_t)}^l$ is an hyperbolic polynomial of degree $l$ and $h(X)$ elliptic of degree $n - l$. Then we denote by $u_i$ and $v_i$ the real and the imaginary part (assumed positive) of complex conjugate zeros of $h(X)$ and by $x_i$ the roots of $g(X)$. Note, that the form does not depend on the order of pairs of conjugate roots. For each cell $\mathcal{P}ol_{(a_1,\ldots,a_t)}^n$ we call the orientation given by $\xi$ the *standard* orientation. We write $\mathcal{P}ol_{(a_1,\ldots,a_t)}^{n,\xi}$ to denote the cell $\mathcal{P}ol_{(a_1,\ldots,a_t)}^n$ oriented by the standard orientation.

Having fixed the above orientations of strata we describe the boundary operator in terms of the incidence coefficients.

**Proposition 2.5.** *Let $(a_1, \ldots, a_t) \preceq (b_1, \ldots, b_{t+1})$ be two compositions in $C_{\leq n}$. Assume $a_i = b_i$, $1 \leq i \leq l$, $a_i = b_{i+1}$, $l+1 \leq i \leq t$, $b_{l+1} = 2$. Let $j_1 \leq l$ be the minimal index such that $b_{j_1} = \cdots = b_{l+1} = 2$ and let $j_2$ be the maximal index such that $b_{l+1} = \cdots = b_{j_2} = 2$. Then:*

(i) $\dim(\mathcal{P}ol_{(a_1,\ldots,a_t)}^n) = \dim(\mathcal{P}ol_{(b_1,\ldots,b_{t+1})}^n) + 1$ *and $\mathcal{P}ol_{(b_1,\ldots,b_{t+1})}^n$ lies in the boundary $\partial \mathcal{P}ol_{(a_1,\ldots,a_t)}^n$.*

(ii) *The coefficient of $\mathcal{P}ol^{n,\xi}_{(b_1,\ldots,b_{t+1})}$ in the expansion of $\partial\mathcal{P}ol^{n,\xi}_{(a_1,\ldots,a_t)}$ equals $0$ if $j_2 - j_1$ is even and $(-1)^{j_1-1}$ if $j_2 - j_1$ is odd.*

**Proof.**
  (i) This follows from Proposition 2.1 and the fact that in $\mathcal{P}ol^{n}_{(b_1,\ldots,b_{t+1})}$ we impose one additional condition on the set of roots.
  (ii) Take a generic polynomial $f(X) \in \mathcal{P}ol^{n}_{(b_1,\ldots,b_{t+1})}$ and a sufficiently small neighborhood $O_f \subset \mathcal{P}ol^{n}_{(b_1,\ldots,b_{t+1})}$. Then $O_f$ lies in the boundary of $j_2 - j_1$ different neighborhoods $O^1_f \subset \mathcal{P}ol^{n}_{(a_1,\ldots,a_t)}, \ldots, O^{j_2-j_1}_f \subset \mathcal{P}ol^{n}_{(a_1,\ldots,a_t)}$. For $1 \leq i \leq j_2 - j_1$ the neighborhood $O^i_f$ is obtained by splitting the double root corresponding to the part $a_i$ in the composition $(a_1, \ldots, a_t)$ into complex conjugate roots for all polynomials from $O_f$. Then the coefficient of $O_f$ in the expansion of $\partial O^i_f$ is $(-1)^{i-1}$. Indeed, $O_f \subseteq \partial O^i_f$ with the orientation

$$dx_i \wedge du_1 \wedge dv_1 \wedge \cdots \wedge du_{(n-\sum a_i)} \wedge dv_{(n-\sum a_i)} \wedge dx_1 \wedge \cdots \wedge \widehat{dx_i} \wedge \cdots \wedge dx_t.$$

Note, that by the definition of the standard orientation the orientation of $O_f \subseteq \partial O^i_f$ does not depend on the choice of the complex root that is split. This orientation differs from the standard orientation by $(-1)^{i-1}$. Taking the sum over all $i$ then implies the assertion. □

Again fixing $(a_1, \ldots, a_t)$ let us denote by $(a_1, \ldots, a_{i-1}, a_i + a_{i+1}, a_{i+2}, \ldots, a_t)$ the sequence of length $t-1$ obtained by merging $a_i$ and $a_{i+1}$. Then the standard stratum $\mathcal{P}ol^{n}_{(a_1,\ldots,a_{i-1},a_i+a_{i+1},a_{i+2},\ldots,a_t)}$ lies in the boundary of $\mathcal{P}ol^{n}_{(a_1,\ldots,a_t)}$ and

$$\dim(\mathcal{P}ol^{n}_{(a_1,\ldots,a_t)}) = \dim \mathcal{P}ol^{n}_{(a_1,\ldots,a_{i-1},a_i+a_{i+1},a_{i+2},\ldots,a_t)} + 1.$$

**Proposition 2.6.** *The coefficient of $\mathcal{P}ol^{n,\xi}_{(a_1,\ldots,a_{i-1},a_i+a_{i+1},a_{i+2},\ldots,a_t)}$ in the expansion of $\partial\mathcal{P}ol^{n,\xi}_{(a_1,\ldots,a_t)}$ equals $(-1)^i$.*

**Proof.** As defined the orientation of strata and differential coincides with the the usual orientation of faces and differential of the boundary of an $n-1$-simplex. □

Now we are in position to state the main result of this section.

**Theorem 2.7.** *Let $\lambda = (\lambda_1 \leq \cdots \leq \lambda_t)$ be a number-partition of a number $l \leq n$ with $n - l$ even.*
  (a) *The homology of $\mathcal{P}ol^{n}_\lambda \cap S^{n-2}$ coincides with homology of the chain complex given by the cells $\mathcal{P}ol^{n,\xi}_{(a_1,\ldots,a_t)}$, for $(a_1, \ldots, a_t) \in C_{\lambda,\leq n}$, as generators and with the differential given in Propositions 2.5 and 2.6.*
  (b) *The homology of $\widehat{\mathcal{P}ol^{n}_\lambda}$ equals the homology of the chain complex given by the cells $\mathcal{P}ol^{n,\xi}_{(a_1,\ldots,a_t)}$, for $(a_1, \ldots, a_t) \in C_{\lambda,\leq n}$, as generators and with the differential given in Propositions 2.5 and 2.6 shifted twice to the right.*

**Proof.** The assertion follows immediately from Proposition 2.5, Proposition 2.6 and Corollary 2.4. □

The next figure illustrates the chain complex for $\lambda = (4)$ in the space of monic polynomials of degree $n = 16$.

**Figure 1 :** The poset $C_{(8), \leq 16}$

**2.4. Stabilization of standard strata.** Let us briefly discuss the properties of $\mathcal{P}ol_\lambda^n$ for a given $\lambda$ and varying $n$. Obviously, there is an inclusion of closures $\mathcal{P}ol^n \subset \mathcal{P}ol^m$ for $n < m$ and $m - n$ even by multiplying each $f(X) \in \mathcal{P}ol^n$ by $(x^2 + 1)^{\frac{(m-n)}{2}}$. By considering the above complex for calculation of $H_*(\widehat{\mathcal{P}ol_\lambda^n})$ one sees by Alexander duality that the cohomology of the complement $H^*(\mathcal{P}ol^n \setminus \mathcal{P}ol_\lambda^n)$ stabilizes. This stabilization is apparently induced by the above inclusion and for $2\sum_{i=1}^{t} \lambda_i - t > 2$ one gets stabilization of the homotopy groups as well. Probably the most interesting detail of this stabilization is that one gets a universal infinite purely combinatorial cell (or chain) complex for the calculation of $H^*(\mathcal{P}ol^\infty \setminus \mathcal{P}ol_\lambda^\infty)$. The cells in this complex are index by compositions in $C_{\lambda,<\infty} = \bigcup_{l<\infty} C_{\lambda,\leq l}$ (the case $\lambda = (3,2)$ is illustrated in Figure 2).

**Figure 2 :** The poset $C_{(2,3),<\infty}$

The situation with the stabilization resembles the one considered in [2] and [25]. One can speculate if the above infinite universal complex associated to each number-partition has another interpretations.

3. HYPERBOLIC POLYNOMIALS, PARTITIONS AND COMPOSITIONS

In this section we will study the special case of the standard strata $\mathcal{H}yp_\lambda$ consisting of hyperbolic polynomials. Before we begin with the study of the spaces $\mathcal{H}yp_\lambda$ we have to introduce some notions and facts about the combinatorics of compositions and partitions.

3.1. **Combinatorics of compositions.** Let $\delta^n$ be the boundary complex of the $n$-simplex (i.e., the simplicial complex spanned by the $(n-1)$-dimensional faces of the $n$-simplex). In particular, $\delta^n \cong S^{n-1}$. We fix an arbitrary bijection $\varphi(\cdot)$ of the $(n+1)$ maximal simplices of $\delta^n$ and the collection of set-partitions

$$\left\{ |1|\cdots|i\ i+1|\cdots|n+2| \;\middle|\; 1 \leq i \leq n+1 \right\}.$$

Thus we have labeled each maximal simplex $\tau$ of $\delta^n$ by a set-partition $\varphi(\tau)$.

Set-partitions are ordered by refinement (i.e., a set-partition $\tau = B_1|\cdots|B_t$ of $[n]$ is less than a set-partition $\sigma = C_1|\cdots|C_s$ of $[n]$ if for all $1 \leq i \leq t$ there is an index $1 \leq j \leq s$ such that $B_i \subseteq C_j$). We write $\tau \preceq \sigma$ if $\tau$ is less than $\sigma$. It is a well known fact that the partially ordered set $\Pi_n$ of all set-partitions of $[n]$ ordered by "$\preceq$" is a lattice (i.e., the supremum operation – denoted by "$\vee$" – and the infimum operation – denote by "$\wedge$" – are well defined).

The labeling $\varphi$ of the $(n-1)$-dimensional faces induces a labeling of all simplices of $\delta^n$ by assigning to a simplex $\sigma$ the join $\bigvee_{\substack{\sigma \subseteq \tau \\ dim(\tau) = n-1}} \varphi(\tau)$ in the lattice $\Pi_{n+2}$ of set-partitions of the set $[n+2]$.

**Lemma 3.1.**
   (i) *Let $\pi$ be a set-partition of $[n]$. Then the following are equivalent:*
      (a) *After possibly reindexing the blocks we have $\pi = B_1|\cdots|B_t$ and $max(B_i) < min(B_{i+1})$, $1 \leq i \leq t-1$.*
      (b) $\pi = \bigvee_{\sigma \in \Gamma} \varphi(\sigma)$ *for some set $\Gamma$ of $(n-3)$-simplices in $\delta^{n-2}$.*
   (ii) *Let $\Gamma_1$ and $\Gamma_2$ be two subsets of the set of $(n-3)$-simplices in $\delta^{n-2}$. Then the following are equivalent:*
      (a) $\Gamma_1 = \Gamma_2$.
      (b) $\bigcap_{\sigma \in \Gamma_1} \sigma = \bigcap_{\sigma \in \Gamma_2} \sigma$.
      (c) $\bigvee_{\sigma \in \Gamma_1} \varphi(\sigma) = \bigvee_{\sigma \in \Gamma_2} \varphi(\sigma)$.

**Proof.**
   (i) If $\pi = B_1|\cdots|B_t$ then define $\Gamma$ by $\sigma \in \Gamma$ if and only if there exists an $1 \leq l \leq n$ and $i, i+1 \in B_l$ such that $\varphi(\sigma) = |1|\cdots|i\ i+1|\cdots|n|$. Then $\bigvee_{\sigma \in \Gamma} \varphi(\sigma) = \pi$. The converse is obvious from the definitions.
   (ii)
   (a) $\Leftrightarrow$ (b) This is easily seen from the the fact that the face lattice of the simplex is isomorphic to the Boolean lattice (see also Remark 3.2 below).
   (a) $\Rightarrow$ (c) This is obvious.
   (c) $\Rightarrow$ (a) If $B$ is a block of $\bigvee_{\sigma \in \Gamma_1} \varphi(\sigma)$. Then for $i < j \in B$ there must be simplices $\sigma_i, \ldots, \sigma_{j-1}$ such that $\varphi(\sigma_l) = |1|\cdots|l\ l+1|\cdots|n|$ for $i \leq l \leq j-1$. But this implies that $\bigvee_{\sigma \in \Gamma_1} \varphi(\sigma)$ determines $\Gamma_1$. □

We call a set-partition a *composition* if it satisfies one of the two equivalent conditions of Lemma 3.1 (ii). This definition does not conflict with the definition of composition given in the introduction. There is a natural bijection sending a set-partition of $[n]$ which is a composition to a composition of the integer $n$ in the usual sense. Namely, let $B_1, \ldots, B_n$ be the reordering of the blocks of a composition $\pi$ such that $max(B_i) < min(B_{i+1})$, then we map $\pi$ to the composition $n = |B_1| + \cdots + |B_t|$. By $C_n$ we denote the partially ordered set – poset for short – of compositions of $n$ ordered by refinement (i.e., by the ordering induced from the usual ordering "$\preceq$" of the lattice $\Pi_n$). Let $\lambda$ be an number-partition of $n$. We denote by $C_\lambda$ the join-semilattice of $C_n$ spanned by the compositions $\pi$ of type $\lambda$ in $C_n$. Of course $C_{n+2} = C_\lambda$ for $\lambda = (2, 1^{n-2})$.

Now we need to introduce some basic material from the theory of partially ordered sets and its order complexes. There is a canonical way to associate to a poset $P$ a simplicial complex $\Delta(P)$. The complex $\Delta(P)$ is the complex of all chains $x_0 < \cdots < x_n$ in $P$ and is called the *order complex* of $P$.

If $\Delta$ is a simplicial complex. Then $\Delta$ may be regarded as a poset ordered by inclusion of faces. The order complex $\Delta(\Delta)$ is the *barycentric subdivision* of $\Delta$ and hence $\Delta(\Delta)$ is homeomorphic to $\Delta$. Note that we assume $\emptyset \notin \Delta$ for a simplicial complex $\Delta$.

By $\delta_\lambda$ we denote simplicial complex spanned by the faces $\tau$ of $\delta^{n-1}$ such that $\varphi(\tau)$ is of type $\lambda$.

**Remark 3.2.**
  (i) The Boolean lattice $B_{n-1}$, the lattice $C_{(1^{n-2},2)}$ and the face lattice of $\delta^{n-1} = \delta_{(1^{n-2},2)}$ are isomorphic lattices.
  (ii) The simplicial complex $\delta_\lambda$ is homeomorphic to a sphere for $\lambda = (1, 2^{n-2})$.

More delicate is the situation for $\lambda \neq (1^{n-2}, 2)$. Before we can proceed, we need some basic tools to manipulate posets preserving the homotopy type. By a map $f : P \to Q$ of posets we understand an *order preserving* map (i.e., $x \leq y$ implies $f(x) \leq f(y)$). A map $f : P \to P$ of a poset into itself is called a *closure operator* if $f(x) \geq x$ and $f(f(x)) = f(x)$. The following lemma is a standard tool in topological combinatorics, we refer the reader for references, proofs and further result to the excellent survey by Björner [4].

**Lemma 3.3.** *Let $P$ be a poset and let $f : P \to P$ be a morphism of posets such that $f(x) \geq x$ and $f \circ f = f$. Then $\Delta(P)$ and $\Delta(f(P))$ are homotopy equivalent.*

Now we are position to state the lemma.

**Lemma 3.4.** *The lattice $C_\lambda$ is homotopy equivalent to the simplicial complex $\delta_\lambda$.*

**Proof.** Let us regard $\delta_\lambda$ as a partially ordered set, ordered by reversed inclusion. Then $\Delta(\delta_\lambda) \cong \delta_\lambda$ (recall, $\Delta(\delta_\lambda)$ is the barycentric subdivision of $\delta_\lambda$). For any $\sigma \in \delta_\lambda$ let $f(\sigma)$ be the union of all maximal simplices of $\delta_\lambda$ containing $\sigma$. Then
  • $f(\sigma) \geq \sigma$,

- $\tau \leq \sigma$ implies $f(\tau) \leq f(\sigma)$,
- $f(f(\sigma)) = f(\sigma)$.

Hence $f$ is a closure operator on the partially ordered set $\delta_\lambda$. Then by Lemma 3.3 it follows that $f$ induces a homotopy equivalence onto its image. But the image of $f$ is by Lemma 3.1 isomorphic to $C_\lambda$. □

3.2. **Stratification and (co)-homology.** The following result links the combinatorics of $\delta_\lambda$ and $C_\lambda$ to the topology of the space $\widehat{\mathcal{H}yp^n_\lambda}$.

**Theorem 3.5.** *Let $\lambda$ be a number-partition of $n$. Then:*

(a) *The one-point compactification $\widehat{\mathcal{H}yp^n_\lambda}$ is homeomorphic to the double suspension over $\delta_\lambda$. In particular, the homeomorphism type of $\widehat{\mathcal{H}yp^n_\lambda}$ depends only on the combinatorics of $\lambda$.*
(b) *The one-point compactification $\widehat{\mathcal{H}yp^n_\lambda}$ is homotopic to to the double suspension over $\Delta(C_\lambda)$. In particular, if $C_\lambda \cong C_{\lambda'}$ for two number-partitions $\lambda$ of $n$ and $\lambda'$ of $m$ then $\widehat{\mathcal{H}yp^n_\lambda} \simeq \widehat{\mathcal{H}yp^m_{\lambda'}}$.*

**Proof.** From part (a) the assertion (b) follows immediately by Lemma 3.4. Thus it remains to verify part (a) of the assertion. We recall, that $\mathcal{H}yp^n$ is homeomorphic to a closed simplicial $(n-1)$-cone $\mathbf{R}^{n-1}_+$ times a copy of $\mathbf{R}$. Thus $\mathcal{H}yp^n_\lambda$ is the cone over $\mathcal{H}yp^n_\lambda \cap S^{n-2}_+$ times a copy of $\mathbf{R}$. Here $S^{n-2}_+$ denotes the positive part of the the unit sphere in $\mathbf{R}^{n-1}$ (i.e., $S^{n-2} \cap \mathbf{R}^{n-1}_+$). The one-point compactification of $Cone(\mathcal{H}yp^n_\lambda \cap S^{n-2}_+)$ is the suspension over $\mathcal{H}yp^n_\lambda \cap S^{n-2}_+$. The one-point compactification of $\mathbf{R} \times \widehat{\mathcal{H}yp^n_\lambda \cap S^{n-2}_+}$ is the suspension over $\left(Cone(\mathcal{H}yp^n_\lambda \cap S^{n-2}_+)\right)^{\widehat{}}$. Since

$$\left(\mathbf{R} \times \widehat{\mathcal{H}yp^n_\lambda \cap S^{n-2}_+}\right)^{\widehat{}} \cong \widehat{\mathcal{H}yp^n_\lambda},$$

the assertion follows one we have identified $\mathcal{H}yp^n_\lambda \cap S^{n-2}_+$ with $\delta_\lambda$.

The intersection of the $(n-1)$-cone $\mathbf{R}^{n-1}_+$ with $S^{n-2}$ is an $(n-2)$-simplex. Its interior corresponds to $n$-tuples $x_1 < \cdots < x_n$ of strictly ascending positive real numbers (on $S^{n-2}$). Its maximal faces correspond to tuples $x_1 < \cdots < x_i = x_{i+1} < \cdots < x_n$ with exactly one equality. In particular, they are naturally labeled by compositions $|1| \cdots |i\ i+1| \cdots |n|$. This in turn establishes the desired homeomorphism of $\mathcal{H}yp^n_\lambda \cap S^{n-2}_+$ with $\delta_\lambda$. □

Now, we describe combinatorially the inclusion $C_\lambda \hookrightarrow B_{n-1}$ that is induced by a fixed isomorphism $C_{(1^{n-2},2)} \cong B_{n-1}$. First we identify the composition $n = a_1 + a_2 + \cdots + a_t$ with the subset $A$ of all elements $i$ of $[n-1]$ such that $a_1 + \cdots + a_j \leq i \leq a_1 + \cdots + a_{j+1} - 1$ for some $0 \leq j \leq t-1$ – the sum over the empty index set is treated as 0. Let us denote this map between $C_{(1^{n-2},2)}$ by $\gamma$. Obviously, $\gamma$ is order preserving and bijective. Thus $\gamma$ induces the desired isomorphism between $B_{n-1}$ and $C_{(1^{n-2},2)}$. Let us study the restriction of $\gamma$ to $C_\lambda$ for $\lambda$ different from $(1^{n-2}, 2)$. Let $\lambda$ be a number-partition of $n$. We say that a subset $A$ of $[n-1]$ is of type $\lambda$ if its preimage under $\gamma$ is a composition of type $\lambda$. In particular, the subsets of type $(1^{n-2}, 2)$ are the singleton subsets of $[n-1]$. For $\lambda = (1^{n-k}, k)$ the subsets of $\lambda$ are the $k-1$ element subsets of $[n-1]$ consisting of $k-1$ consecutive numbers. More generally, our considerations so far imply.

**Lemma 3.6.** *$C_\lambda$ is isomorphic to the join-subsemilattice of $B_{n-1}$ that is generated by all subsets of $[n-1]$ of type $\lambda$.*

If $n = a_1 + \cdots + a_t$ is a composition for a number $n$ there may possibly be more than two ways of expressing a number $l \leq n$ as a sum using the parts $a_1, \ldots, a_t$. For example take the composition $4 = 2 + 1 + 1$ then $2 = 1 + 1$ is such an identity that holds within the given composition. An identity $n_1 + \cdots + n_r = m_1 + \cdots + m_s$ is called *primitive partition identity* (see for example Sturmfels' book [23, Chap. 6]), if $n_i \neq m_j$ for $1 \leq i \leq r$ and $1 \leq j \leq s$. We say that a composition $n = a_1 + \cdots + a_t$ is *free* of *resonances* if there is no primitive partition identity $\sum_{i \in I} a_i = \sum_{j \in J} a_j$, for disjoint non-empty subsets $I, J$ of $[t]$. Clearly, the property of being free of resonances only depends on the type of the composition. So one may as well apply "free of resonances" to number-partitions. For the description of the $\Delta(C_\lambda)$ in case $\lambda$ is free of resonances, we use an important class of polytopes. The *permutahedron* $P_t$ is the $(t-1)$-dimensional polytope – for us the polytope is the sphere subdivided by the polytope – spanned in $\mathbf{R}^t$ be the the vectors $(\sigma(1), \ldots, \sigma(t))$ for $\sigma \in S_t$ (note that all these points lie on the hyperplane $x_1 + \cdots + x_t = \frac{t(t-1)}{2}$. Alternatively, the permutahedron allows the following description. Let $\mathcal{A}_t = \{\tilde{H}_{ij} : x_i = x_j \mid 1 \leq i < j \leq t\}$ be the arrangement (i.e., finite set) of reflecting hyperplanes of the Coxeter group $S_t$ in $\mathbf{R}^t$. Let $H : x_1 + \cdots + x_t = 0$ be the hyperplane orthogonal to the line $x_1 = \cdots = x_t$. Then the arrangement triangulates the intersection $S^{t-2}$ of $H$ with the unit sphere in $\mathbf{R}^t$. Let $\Delta_t$ be this triangulation of $S^{t-2}$. Then let $P_t$ the cellular dual of $\Delta_t$ (i.e., the cellular complex constructed by replacing cells of dimension $i$ in $\Delta_t$ by a cell of dimension $t-2-i$ preserving adjacency; see the book by Björner, Las Vergnas, Sturmfels, White & Ziegler [5] for more details on permutahedra). One easily checks that the cellular dual of a polytope again can be realized as a polytope, thus $P_t$ is indeed a polytope. We denote by $L_t$ the poset of faces of $P_t$ ordered by reverse inclusion. Again it is a well known fact that $L_t$ is $S_t$-isomorphic to the poset of cosets of Young-subgroups $\neq S_t$ of $S_t$ ordered by inclusion. In general if $P$ is a poset on which a group $G$ acts as a group of (poset)-automorphisms then we denote by $P/G$ the poset on the set of orbits $x^G = \{x^g \mid g \in G\}$ whose order relation is given by $x^G \leq y^G :\Leftrightarrow \exists\, g \in G : x^g \leq y$.

**Proposition 3.7.** *Let $\lambda = (1^{e_1}, \ldots, n^{e_n})$ be a number-partition of $n$ that is free of resonances. Let $H = S_{e_1} \times \cdots \times S_{e_n}$ be the corresponding Young-subgroup of $S_t$, $t = e_1 + \cdots + e_n$. Then:*

(i) *The poset $C_\lambda$ is isomorphic to the lattice of cells of the quotient $P_t/H$ of the permutahedron $P_t$. In particular, $C_\lambda$ is homotopic to $P_t/H$.*
(ii) *If $0 \le e_i \le 1$ then $C_\lambda$ is $S_t$-isomorphic to the face lattice of the permutahedron $P_t$, $t = e_1 + \cdots + e_n$. The order complex $\Delta(C_\lambda)$ triangulates an $(t-2)$-sphere.*
(iii) *If $e_i \ge 2$ for some $i$ then $P_t/H$ and $\Delta(C_\lambda)$ are contractible.*

**Proof.**

(i) By the second description of the permutahedron given above, the quotient of the permutahedron $P_t$ by a Young-subgroup $S_{b_1} \times \cdots \times S_{b_f}$ is homeomorphic to the closure of the intersection of $P_t$ with a fundamental domain $C$ of the Coxeter group $S_{b_1} \times \cdots \times S_{b_f}$ in $\mathbf{R}^t$. Thus if $F$ is a (closed) face of $P_t$ that intersects the fundamental domain $C$ then the image of $F$ under the projection map $P_t \to P_t/H$ is $F \cap \overline{C}$. Thus the images of the (closed) faces of $P_t$ under the action of $H$ are ordered along the orbits poset of the poset $L_t$ under the action of $H$.

If $X$ is a CW-complex that is stratified into closed cells $C_p$, $p \in P$ for a poset $P$, such that:
(a) The set $\{C_p \mid p \in P\}$ is closed under taking intersections.
(b) The order relation in $P$ is the reversed inclusion relation on $\{C_p \mid p \in P\}$.
Then $X \simeq \Delta(P)$. Since (a) and (b) are satisfied for $X = P_t/H$ and $P = L_t/H$ the assertion follows once we have verified that $C_\lambda$ is isomorphic to $L_t/H$.

Obviously, if $a_1 + \cdots + a_t$ is a composition of type $\lambda$ and $\sigma, \tau \in S_t$ are permutations then $a_{\sigma(1)} + \cdots + a_{\sigma(t)}$ is a composition different from $a_{\tau(1)} + \cdots + a_{\tau t}$ if and only if $\sigma\tau^{-1} \notin H$. Thus the minimal elements of $C_\lambda$ and the minimal elements of $L_t/H$ are in bijection. Assume, we are given a composition $b_1 + \cdots + b_s$ in $C_\lambda$. Let $\sigma, \tau$ be two permutations such that $a_{\sigma(1)} + \cdots + a_{\sigma(t)}$ and $a_{\tau(1)} + \cdots + a_{\tau(t)}$ are finer than $b_1 + \cdots + b_s$. Since there are no resonances, it follows that $\sigma$ and $\tau$ must lie in the same $H$ orbit on the cosets of $S_{b_1} \times \cdots \times S_{b_s}$ in $S_t$. Conversely, if $\sigma$ and $\tau$ lie in the same $H$-orbit on the cosets of $S_{b_1} \times \cdots \times S_{b_s}$ in $S_t$ but in different cosets of $H$ in $S_t$ then $a_{\sigma(1)} + \cdots + a_{\sigma(t)}$ and $a_{\tau(1)} + \cdots + a_{\tau(t)}$ are different. This induces a bijection between the minimal elements below $c_1 + \cdots + c_s$ in $C_\lambda$ and the minimal elements below an $H$ orbit on the cosets of $S_{b_1} \times \cdots \times S_{b_s}$ in $S_t$. Since in $C_\lambda$ and $L_t/H$ each element is the supremum of all minimal elements below it, the assertion follows.

(ii) If $H = 1$ then $L_t/H = L_t$ and $P_t/H = P_t$. Thus $P_t/H$ and $\Delta(C_\lambda)$ triangulate a $(t-2)$-sphere.

(iii) If $e_i \geq 2$ for some $e_i$ then the intersection of $P_t$ with the fundamental domain $C$ (i.e., a convex cone) is homeomorphic to a $(t-2)$-ball. Thus $C_\lambda$ is contractible. □

**Theorem 3.8.** *Let $\lambda = (1^{e_1}, \ldots, n^{e_n})$ be a number-partition of $n$ that is free of resonances. Then $\widehat{\mathcal{H}yp}^n_\lambda$ is contractible if $e_i \geq 2$ for some $i$ and homeomorphic to an $(e_1 + \cdots + e_n)$-sphere if $0 \leq e_i \leq 1$, $1 \leq i \leq n$.*

**Proof.** The assertion follows immediately from Proposition 3.7 and Theorem 3.5. Note that the suspension of a $i$-sphere is homeomorphic to an $(i+1)$-sphere. □

We come to explicit results on the homotopy type for special classes of number-partitions $\lambda$. We consider number-partitions $\lambda = (1^{n-k}, k)$, $2 \leq k \leq n$. In this case the subsets of $[n-1]$ of type $\lambda$ are intervals themselves. Thus $C_\lambda$ is isomorphic to the sub-poset of $B_{n-1}$ that is generated by the intervals $[(k-1)\cdot t+1, (k-1)\cdot(t+1)]$ for $0 \leq t \leq \lfloor \frac{n}{(k-1)} \rfloor$. In general, a join-sublattice $L$ of $B_n$ is called *interval generated* if there is a set $I$ of intervals $I = \{[a_i, b_i] \mid i = 1, \ldots, t\}$ such that $L$ is generated by the elements of $I$. Interval generated sublattices of $B_n$ have been introduced and studied by Greene [15]. Recently the method of non-pure shellability, a concept due to by Björner & Wachs [7], [6], allowed to compute in [6] the homotopy type of the lattice of intervals generated by all $(k-1)$-element subsets of $[n-1]$.

**Proposition 3.9.** [6, Theorem 8.2, Corollary 8.4] *Let $\lambda = (1^{n-k}, k)$. Then,*
$$\Delta(C_{(1^{n-k},k)}) \simeq \begin{cases} S^{2(n-1)/k-2}, & \text{if} \quad n \equiv 1 \pmod{k}, \\ S^{2n/k-3}, & \text{if} \quad n \equiv 0 \pmod{k}, \\ \text{pt.}, & \text{otherwise.} \end{cases}$$

**Proof.** By Lemma 3.6 the poset $C_{(1^{n-k},k)}$ is isomorphic to the join-subsemilattice of $B_{n-1}$ generated by all intervals of length $k-1$. Now by the result of Björner & Wachs [6, Theorem 8.2, Corollary 8.4] the assertion follows. □

**Corollary 3.10.** *Let $\lambda = (1^{n-k}, k)$. Then*
$$\widehat{\mathcal{H}yp}^n_\lambda \simeq \begin{cases} S^{2(n-1)/k}, & \text{if} \quad n \equiv 1 \pmod{k}, \\ S^{2n/k-1}, & \text{if} \quad n \equiv 0 \pmod{k}, \\ \text{pt.}, & \text{otherwise.} \end{cases}$$

**Proof.** The assertion follows immediately from Proposition 3.9 together with Theorem 3.5. □

Finally, we would like to provide the results of some computer computations using a program package for Mathematica© written by Vic Reiner. For the parameters $\lambda = (1^{n-4}, 2^2)$, $n \leq 9$, $\lambda = (1^{n-6}, 3^2)$ for $n \leq 11$ all homology of $\Delta(C_\lambda)$ vanishes. For $\lambda = (a, b, c, d)$, $a < b < c < d$, with $a + b = c$, $b + c = d$ as the only primitive partition identities, homology is free of rank 3 and concentrated in dimension 2. For $\lambda = (a, b, c, d)$, $a < b < c < d$, with $a + b + c = d$ as the only primitive partition identity, homology is free of rank 1 in dimensions 1 and 2.

3.3. **Iterated compositions.** One may regard the consideration of hyperbolic polynomials as the 1-dimensional case of a general situation. To understand the general philosophy, one has to regard $\mathcal{H}yp^n$ as the space of unordered point configurations (i.e., configurations of roots) on the real line $\mathbf{R}^1$. The essential combinatorial ingredient for the 1-dimensional case – compositions – will prevail to be important in higher dimensions too. Most of what follows now is a translation of Section I.3.3 of Vassiliev's book [25] into the language of posets and compositions. Let $d \geq 1$ be some number. The orbit space $(\mathbf{R}^d)^n/S_n$ admits the following cell decomposition. An *iterated-composition* of degree $d$ of the number $n$ is a $d$-tuple $(\pi_1 \preceq \cdots \preceq \pi_d)$ of compositions of $n$. We order the iterated compositions of $n$ of degree $d$ by the order induced by the lexicographic order on $d$-tuples of compositions. By $C_n^d$ we denote the poset of iterated compositions of $n$ of degree $d$.

**Proposition 3.11.** *The poset $C_n^d$ of iterated compositions of $n$ of degree $d$ is isomorphic to a direct product of $n-1$ chains of cardinality $d+1$.*

**Proof.** We already know that the poset of compositions $C_n^1 = C_n$ is isomorphic to the Boolean lattice $B_{n-1}$. Thus the assertion follows for the case $d = 1$. Let $\pi_i$, $i = 1, \ldots, n$ be the atoms of the lattice of compositions of $[n]$. Let $\hat{0}$ be the least element of the lattice of compositions of $n$. Then we define chains $c_i : x_{i,0} < \cdots < x_{i,d}$ by setting
$$x_{i,j} = (\underbrace{\hat{0} \preceq \cdots \preceq \hat{0}}_{d-j \text{ times}} \preceq \underbrace{\pi_i \preceq \cdots \preceq \pi_i}_{j \text{ times}}).$$
Now one checks that $C_{n+1}^d \cong c_1 \times \cdots \times c_{n-1}$. □

To each iterated composition $\pi$ of $n$ of degree $d$ one associates (see Vassiliev's book [25, II,3.3]) a cell $c(\pi)$ in $(\mathbf{R}^d)^n/S_n$. Namely, let $\pi = (\pi_1 \preceq \cdots \preceq \pi_d)$ be an iterated composition of $n$ of degree $d$. The element $\underline{\mathbf{x}} = (x_1, \ldots, x_d) \in (\mathbf{R}^d)^n/S_n$ is an element of $c(\pi)$ if and only if

(i) for every pair $i, j$ of elements which lie in the same block of $\pi_s$ for some $1 \leq s \leq d$ the coordinates $x_s^i$ and $x_s^j$ satisfy $x_s^i = x_s^j$,
(ii) for every pair $i < j$ of elements which do not lie in the same block of $\pi_s$, but lie in the same block of $\pi_{s-1}$ (here we set $\pi_0 = |1| \cdots |n|$), for some $1 \leq s \leq d$ the coordinates $x_s^i$ and $x_s^j$ satisfy $x_s^i \leq x_s^j$.

**Lemma 3.12.** [25, Lemma 3.3.1, 3.3.2] *Let $\pi = (\pi_1 \preceq \cdots \preceq \pi_d)$ be an iterated composition of $n$ of degree $d$. Assume $\pi_i$ consists of $t_i$ blocks, $1 \leq i \leq d$. Then $c(\pi)$ is a cell of dimension $nd - (n - t_1) - \cdots - (n - t_d)$.*

Again one is able to identify within $C_n^d$ certain strata that correspond to various degenerations in $(\mathbf{R}^d)^n/S_n$. Let $\lambda$ be a number-partition of $n$ then set $C_\lambda^d$ to be the join-sublattice of $C_{n+1}^d$ that is generated by iterated compositions $(\pi \preceq \cdots \preceq \pi)$ for which $\pi$ is of type $\lambda$. The corresponding cells stratify a $d$-dimensional analog of $\mathcal{H}yp^n$. But for general $d$ the situation is more delicate than in the case $d = 1$. The cells $c(\pi)$ for $\pi \in C_\lambda^d$ are attached to each other in a more subtle than in case $d = 1$. Still it is possible to compute the boundary map (see Fuchs [14], Cohen [10], Vainshtain [24] and Vassiliev [25]).

## 4. Concluding remarks and problems

4.1. **Hyperbolic polynomials & resonance arrangements.** By Theorem 3.5 (a) and (b) the topology of $\widehat{\mathcal{H}yp_\lambda^n}$ is determined by the "combinatorics" of $C_\lambda$. On the other hand it is clear that the lattice structure of $C_\lambda$ is determined by the "resonances" of $\lambda$. Note that any lattice identity in $C_\lambda$ corresponds to some resonance identity of $\lambda$ and vice versa. Therefore, it is reasonable to ask for the "set of resonances" valid for a number-partition $\lambda$ or equivalently for a composition $\pi$ of type $\lambda$. This leads necessarily to the study of a certain arrangement of hyperplanes. Assume $\lambda = (1^{e_1}, \ldots, n^{e_n})$, $t = \sum_{i=1}^{n} e_i$. Then we consider a composition $(a_1, \ldots, a_t)$ of type $\lambda$ as a point in $\mathbf{R}^t$. Let $\mathcal{A}_t$ be the arrangement (i.e., finite set) of hyperplanes

$$\Big\{H_{J,K} : \sum_{i \in J} x_i = \sum_{i \in K} x_i \mid \emptyset \neq J, K \subseteq [t], J \cap K = \emptyset\Big\}.$$

For a composition $\pi = (a_1, \ldots, a_t)$ let $\mathcal{A}_\pi$ be the set of hyperplanes in $\mathcal{A}_t$ that contain the point $(a_1, \ldots, a_t)$. Note, that $\mathcal{A}_\pi$ is empty if and only if $a_i \neq a_j$ for $i \neq j$ and $(a_1, \ldots, a_t)$ is free of resonances. Now the combinatorics of $C_\lambda$ is determined by the $S_t$ orbit $\mathcal{A}_\pi$ or equivalently $\bigcap_{H \in \mathcal{A}_\pi} H$ in the positive orthant (i.e., $x_i \geq 0$, $1 \leq i \leq t$), where $S_t$ acts on $\mathbf{R}^t$ by permuting coordinates. Therefore, it appears as a challenging problem to count the number of different $S_t$ orbits. The problem of counting the number of connected components of the complement $\mathbf{R}^t \setminus \bigcup_{H \in \mathcal{A}_t} H$ amounts to the enumeration of the number of different classes of resonance free compositions with no multiple parts. This problem has evoked interest in statistics (see Fine & Gill [13]) and electrical engineering. In general, arrangements that contain only hyperlanes defiened by linear form with 0, 1 and $-1$ as coefficients are were studied by Edelman & Reiner [12].

Another challenging problem that arises here, is the question : "How general are the order complexes $\Delta(C_\lambda)$ for various $\lambda$ ?" The results presented and the examples supplied provide only evidence for "very nice" simplicial complexes (i.e., $\Delta(C_\lambda)$ has the homotopy type of a sphere). So one may as : "Is it true that in general $\Delta(C_\lambda)$ is homotopic a wedge of spheres ? " Or "Is it true that in general the homology of $\Delta(C_\lambda)$ is free ? " See the examples at the end of Section 2.2. On the other hand the number-partitions $(1^{n-k}, k)$ allow only a single and very basic resonance. Namely, $k = \underbrace{1 + \cdots + 1}_{k \text{ times}}$. Thus the topology of $\Delta(C_\lambda)$ may be much more complicated in general. We mention that for $\lambda = (k^r, 1^q)$, $q \geq k \geq 3$, $r \geq 2$ the method of shellability used by Björner & Wachs [6] to determine the homotopy type fails. This follows by a similar argument used by Kozlov [19] for analogously defined partition lattices. So conversely, it may even be true that the complexes $\Delta(C_\lambda)$ are universal in the homotopy category of finite CW-complexes (i.e., for any finite CW-complex $X$ there is a number-partition $\lambda$ such that $X \simeq \Delta(C_\lambda)$). We regard this as a very interesting, albeit difficult, problem. A combinatorial reformulation and generalization of this problem can be found in [22], where "weight"-functions on lattices are considered.

4.2. **Quotients of permutahedra.** The procedure of passing from $P_t$ to $P_t/H$ for a Young-subgroup $H$ of $S_t$ applied in Proposition 3.7 is reminiscent of the construction of the braid space as a quotient of the permutahedron (see Carlsson & Milgram [9]). For the other Artin braid groups this was generalized by Salvetti [21]. Let us briefly sketch and compare the two constructions in order to exhibit the significant difference between them. The construction described by Carlsson & Milgram and Salvetti is a quotient construction under the full symmetric group (resp., Coxeter group). Thus we start describing the cellular complex that is constructed in both cases for the quotient operation under the "action" of $S_t$.

(i) The quotient space of the permutahedron $P_t$ modulo the action of $S_t$ has exactly $\binom{t-1}{i}$ cells of dimension $i$, $0 \leq i \leq t-2$. Namely, by Proposition 3.7 the cellular decomposition of $P_t/S_t$ is given by the cells of a simplex. Note that $L_t/S_t$ is isomorphic to the Boolean lattice on $t-1$ elements with the top element removed.

(ii) The second complex, that arises by "modding" out the $S_t$-action on the permutahedron, consists as well of $\binom{t-1}{i}$ cells of dimension $i$, $0 \leq i \leq t - 2$. But here attachment and the boundary maps are much more complicated. Let $v$ be a vertex of the permutahedron. For each face $F$ of the permutahedron let $\sigma_F$ be the permutation of minimal length – length as a word in Coxeter generators – in $S_t$ such that $F^{\sigma_F}$ intersects the fundamental chamber of the $S_t$ action on $\mathbf{R}^t$ – here we use the second description of the permutahedron. If two faces $F$, $F'$ lie in the same $S_t$ orbit, then they are be identified via $\sigma_F \sigma_{F'}^{-1}$. The resulting space will be denoted by $Braid_t$. It has be shown (see Carlsson & Milgram [9]) that $Braid_t$ is the $K(\pi, 1)$ space for the braid group $\pi = Br_t$ on $t$ strings. For the other Coxeter groups the result is due to Salvetti [21].

In Figure 3 we perform the two constructions for $t = 3$. Obviously, the two constructions do not yield homotopy equivalent spaces.

Classifying space of the braid group $Br_1 \cong \mathbf{Z}$       Permutahedron $P_3$       Quotient space $P_3/S_3$

**Figure 3 :** Quotient constructions on the permutahedron $P_3$

4.3. **Complex polynomials.** One can consider similarly defined stratification by multiplicities of roots of the space of monic polynomials of degree $n$ with complex coefficients. Note that these strata are enumerated by the set of all number-partitions of $n$ and the only stratum of top dimension coincides with the standard braid space. The strata corresponding to a number-partition whose Young diagram is of the hook shape were studied in the work of Arnol'd [1] and Vassiliev [25]. In particular, Vassiliev has shown that the homotopy type of the 1-point compactification of such a stratum is stably homotopically equivalent to a wedge of vector bundles over the braid spaces coming from the sign representation in the symmetric group. These spaces are not so simple to approach and combinatorial stratifications may be more difficult to handle. Nevertheless, teh first author has observed that the homotopy type of each such stratum depends only on the combinatorics of the corresponding number-partition. More generally, the first author has checked that the homotopy type of $\widehat{\Sigma_{\mu,\lambda}}$ depends only on the combinatorics of the number-partitions $\mu$ and $\lambda$. But in these cases homology and cohomology computations become much more delicate and torsion appears all over the place (see the work of Fuchs [14], Cohen [10], Vainshtain [24] for the case $\lambda = (1^{n-2}, 2)$). It will be very interesting to find the analog of Vassiliev's result for the case of any number-partition.

4.4. **Stable complexes.** The problem of calculation of the stable cohomology of the complement $\mathcal{P}ol^\infty \setminus \mathcal{P}ol_\lambda^\infty$ mentioned in Section 2 is closely related to the one considered in [2] and [25]. One can hope that this question is solvable for any number-partition $\lambda$. In particular, if each entry of $\lambda$ is greater or equal 3 then the technique of Vassiliev's resolving spaces [25] is applicable and one can apparently get a (simple?) answer for the additive structure of $H^*(\mathcal{P}ol^\infty \setminus \mathcal{P}ol_\lambda^\infty)$.

## References

[1] V.I. Arnol'd. Topological invariants of algebraic functions. *Trans. Moscow Math. Soc.*, 21:30–52, 1970.
[2] V.I. Arnol'd. Spaces of functions with moderate singularities. *Functional Anal. Appl.*, 23(3):169–177, 1989.
[3] A. Björner. Subspace arrangements. In A. Joseph, editor, *Proceedings of the 1st European Congress of Mathematics (Paris 1992)*, volume 1, pages 321–371, Boston, 1995. Birkhäuser.
[4] A. Björner. Topological methods. In R. Graham, M. Grötschel, and L. Lovász, editors, *Handbook of Combinatorics*, pages 1819–1872. North-Holland, Amsterdam, 1995.
[5] A. Björner, M. Las Vergnas, B. Sturmfels, N. White, and G.M. Ziegler. *Oriented Matroids*. Encyclopedia of Mathematics. Cambridge University Press, Cambridge, 1993.
[6] A. Björner and M.L. Wachs. Shellable non-pure complexes and posets, II. *Trans. Amer. Math. Soc.*, 1995. (to appear).
[7] A. Björner and M.L. Wachs. Shellable non-pure complexes and posets, I. *Trans. Amer. Math. Soc.*, 348(4):1299–1327, 1996.
[8] R.W. Brockett. Some geometric questions in the theory of linear systems. *IEEE Trans. Automat. Control*, AC-21:449–455, 1976.
[9] G. Carlsson and R.J. Milgram. Stable homotopy and iterated loop spaces. In I.M. James, editor, *Handbook of Algebraic Topology*, pages 505–583. North Holland, Amsterdam, 1995.
[10] F.R. Cohen. The homology of $C_{n+1}$-spaces, $n \geq 0$. In F.R. Cohen, T.J. Lada, and J.P. May, editors, *The homology of iterated loop spaces*, volume 533 of *Lecture Notes in Math.*, pages 207–353. Springer, Heidelberg, 1976.
[11] J. Dedieu. Obreschoff's theorem revisited: what convex sets are contained in the set of hyperbolic polynomials. *J. Pure Appl. Algebra*, 81(3):269–278, 1992.
[12] P. Edelman and V. Reiner. Free hyperplane arrangements between $A_n$ and $B_n$. *Math. Z.*, 215:347–365, 1994.
[13] T. Fine and J. Gill. The enumeration of comparative probability relations. *Ann. Probab.*, 4:667–673, 1976.
[14] D. Fuchs. Cohomologies of the group cos mod 2. *Functional Anal. Appl.*, 4:143–151, 1970.
[15] C. Greene. A class of lattices with Möbius function $\pm 1$, 0. *European J. Combin.*, pages 225–240, 1988.
[16] U. Helmke. A compactification of the space of rational transfer functions by singular systems. *J. Math. Systems*, 3(4):459–472, 1993.
[17] U. Helmke and C.I. Byrnes. The cohomology of the moduli space of controllable linear systems. *Acta Appl. Math.*, 28:161–188, 1992.
[18] L. Hörmander. *The Analysis of Linear Partial Differential Operators*, volume 256 of *Grundlehren der Math. Wissenschaften*. Springer, Heidelberg, 2 edition, 1990.
[19] D. Kozlov. On generalized lexicographic shellability. (Preprint 1995).
[20] I. Meguerditchian. A theorem on escape from the space of hyperbolic polynomials. *Math. Z.*, 211:350–359, 1992.
[21] M. Salvetti. The homotopy type of the Artin groups. *Math. Res. Lett.*, 1:565–577, 1995.
[22] B. Shapiro. Valuations on posets and associated simplicial complexes. University of Stockholm, Preprint 1996.

[23] B. Sturmfels. *Gröbner Bases and Convex Polytopes*, volume 8 of *University Lecture Series*. Amer. Math. Soc., Providence, RI, 1995.
[24] F.V. Vainshtain. Cohomology of the braid groups. *Functional Anal. Appl.*, 12:135–137, 1978.
[25] V.A. Vassiliev. *Complements of Discriminants of Smooth Maps : Topology and Applications*, volume 98 of *Transl. of Math. Monographs*. Amer. Math. Soc., Providence, RI, 1994. Revised Edition.

DEPARTMENT OF MATHEMATICS, STOCKHOLM UNIVERSITY, S-106 91 STOCKHOLM, SWEDEN
*E-mail address*: shapiro@matematik.su.se

INSTITUTE FOR EXPERIMENTAL MATHEMATICS, ELLERNSTR. 29, 45326 ESSEN, GERMANY
*E-mail address*: welker@exp-math.uni-essen.de